\documentclass{ws-aa}

\usepackage{amssymb,wasysym}
\newcommand\R{{\mathbb R}}

\newcommand\Dis{\displaystyle}
\newcommand\fa{\frac}
\newcommand\veps{\varepsilon}

\begin{document}


%
%

\title{H-measures and system of Maxwell's \\}

\author{Hassan TAHA}
\address{Universit\'e d'Orl\'eans \\
MAPMO-UMR 6628, BP 6759\\
45067 ORLEANS CEDEX-France\\ 
hassan.taha@labomath.univ-orleans.fr}
\maketitle

\begin{abstract}
We are interested in the homogenization of energy 
like quantities in electromagnetism.  We prove a general propagation Theorem 
for H-measures associated to Maxwell's system, 
in the full space $\Omega =\R^{3}$, without boundary conditions.
We shall distinguish between two cases: constant coefficient case, and non  coefficient-scalar
case. In the two cases we give the behaviour of the H-measures associated to this 
system. 
\end{abstract}

\keywords{Electromagnetism, homogenization of energy, H-measures, Maxwell's system.}

\ccode{Mathematics Subject Classification 2000: 35BXX, 35B27}

\section{Introduction}
Herein, we are interested in the homogenization of energy 
like quantities in electromagnetism, and more particularly in Maxwell's 
equations, without  boundary conditions. We use the notion of 
H-measures, introduced by G\'erard and Tartar \cite{Gér1}, \cite{tar}. We 
prove a general propagation Theorem for H-measures associated to 
Maxwell's system. This result, combined with the localisation property, is 
then used to obtain more precise results on the behaviour of H-measures 
associated to this system.\par 
\noindent As known, an H-measure is a (possibly matrix of) Radon measures on the 
product space $\Omega\times S^{\;n-1}$, where $\Omega \subseteq \R^{n}$ is 
an open domain and $S^{\;n-1}$ is the unit sphere in $\R^{n}$. In order to apply 
Fourier transform, functions defined on the whole of $\R^{n}$ should be 
considered and this can be achieved by extending them by zero outside 
the domain. For this reason, we consider Maxwell's system in the full space 
$\R^{3}$, which means without boundary conditions. Let us mention that similar works already exist, see in particuliar 
\cite{ant}. However, in \cite{ant}, computations are far from being complete.\par
\noindent If one is interested in coupling with boundary conditions, the usual pseudodifferential calculus behind 
the notion of H-measures is not sufficient, and one should use much more technical tools. In the case of semi-classical 
measures, this is now well known, see for instance \cite{Gér2}, \cite{Géreric}, \cite{lmul}.\par
\noindent In the context of H-measures, and in particuliar without a typical scale, then one should use tools similar to 
those developped in recents works, see for instance \cite{alexa}. However, the results presented here will be important for the full Maxwell problem, with suitable 
boundary conditions.\par
\noindent Let $\Omega$ be an open set of $\R^{\;3}$. We 
will consider Maxwell's system in the material $\Omega$ with electric 
permeability $\ddot{\bf\epsilon} $, conductivity 
$\ddot{\bf\sigma} $ and magnetic susceptibility $\ddot{\bf\eta} $ given by
\begin{equation}\label{maxwel}
\left\{
\begin{array}{ccccccc}
i) &\partial_{t}D^{\Dis\varepsilon} (x,t) +J^{\Dis\varepsilon} (x,t) &= &\mbox{rot} 
H^{\Dis\varepsilon} (x,t) +F^{\Dis\varepsilon},&\\[0.3cm]     
ii) &\partial_{t} B^{\Dis\varepsilon} (x,t) &= &-\mbox{rot} E^{\Dis\varepsilon} 
(x,t) +G^{\Dis\varepsilon} (x,t),&\\[0.3cm]     
iii) &\mbox{div} B^{\Dis\varepsilon} (x,t)& =&0,&\\[0.3cm] 
iv) &\mbox{div} D^{\Dis\varepsilon} (x,t)& =&\varrho^{\Dis\varepsilon} (x,t),&   
\end{array}
\right.
\end{equation}

\noindent where $x\in\Omega$ and $t\in(0,T)$, 
$ E^{\Dis\veps}$, $H^{\Dis\veps}$, $D^{\Dis\veps}$, $J^{\Dis\veps}$ and 
$B^{\Dis\veps}$ are the electric, magnetic, induced electric, current density 
and induced magnetic fields, respectively.\par
\noindent Morever, $ \rho^{\;\Dis\veps}$ (the charge density), $ F^{\;\Dis\veps}$ and 
$ G^{\;\Dis\veps}$ are given, and we have the three constitutive relations  
\begin{equation}\label{lois-const}
\left\{
\begin{array}{ccccccc}    
1) &\Dis D^{\Dis\veps}(x,t)&=& \ddot{\epsilon}(x) E^{\Dis\veps}(x,t),&\\[0.3cm] 
2) &\Dis J^{\Dis\veps}(x,t)&=&\ddot{\sigma}(x)E^{\Dis\veps}(x,t),&\\[0.3cm] 
3) &\Dis B^{\Dis\veps}(x,t)&=&\ddot{\eta}(x)H^{\Dis\veps}(x,t),&\\[0.3cm]      
\end{array}
\right.
\end{equation}\noindent where $\ddot{\epsilon}$, $\ddot{\sigma}$ and $\ddot{\eta}$ are $3\times 3$ matrix valued 
functions and $\varepsilon$ is a typical lenght going to $0$.\par 
\noindent We shall consider this system in the full space $\Omega =\R^{3}$, without boundary 
conditions. Since we are not taking into account the initial data, we will also assume that the time variable $t$ belongs to $\R$.\par
\noindent We shall use the notion of H-measure to compute for instance energy quantities 
in the following cases :\\
\noindent {\bf i) Constant coefficient case}: here, we assume that the 
electric permittivity $\ddot{\Dis\epsilon}$, conductivity $\ddot{\sigma}$ and 
magnetic susceptibility $\ddot{\eta}$  are $3\times 3$ identity matrices, i.e. 
\begin{equation}\label{hyp1}
\ddot{\Dis\epsilon}=\ddot{\sigma}=\ddot{\eta}=(Id)_{3\times 3}\;.
\end{equation}
\noindent {\bf ii) Non constant coefficient-scalar case}: in this case, we 
consider that the matrix $\ddot{\Dis\epsilon}$, $\ddot{\sigma}$,  $\ddot{\eta}$ are scalar
$3\times 3$ matrix valued smooth functions, i.e. 
\begin{equation}\label{hyp2}
\ddot{\Dis\epsilon}= \Dis\epsilon (Id)_{3\times 3}, \hspace{0.2cm} \ddot{\sigma}=\sigma (Id)_{3\times 3}
,\hspace{0.2cm} \ddot{\eta}=\eta (Id)_{3\times 3}
\end{equation}
\noindent where $\Dis\epsilon$, $\sigma$,  $\eta$ are smooth functions, given in $C^{1}_b(\R^{3})$, bounded from below.\par
\noindent We will also assume that
\begin{equation}\label{hyp-conv-donnees}
f^{\Dis\veps}=(F^{\Dis\veps}, G^{\Dis\veps})^{t},\; \varrho^{\Dis\veps} \rightharpoonup (0,0,0) \mbox{ in } [L^2(\R\times \R^3)^6]\times L^2(\R\times \R^3) \mbox{ weakly }
\end{equation}
and
\begin{equation}\label{hyp-conv-u}
u^{\Dis\veps} \equiv (E^{\Dis\veps} ,H^{\Dis\veps})^t \rightharpoonup 0 \mbox{ in } L^2 (\R\times \R^3)^6 \mbox{ weakly.}
\end{equation}
\noindent  After some prerequisites on H-measures, see \cite{Gér2} and \cite{tar}, presented in Section I, we use this notion in Section II to prove 
\noindent\begin{theorem}
\label{theoreme1}{\bf Constant coefficient case}\\
\noindent Assume (\ref{lois-const}), (\ref{hyp-conv-donnees}), (\ref{hyp-conv-u}) and (\ref{hyp1}). Then, up to a suitable extraction, the H-measure $\scorpio = \scorpio (t,x,\zeta )$, $\zeta =(\zeta_0 ,\zeta ')$, $\zeta ' =(\zeta_1 , \zeta_2 , \zeta_3 )$, associated to $(u^{\Dis\veps})$, can be expressed as follows  
\begin{equation}\label{1.8}
\scorpio =\left(
           \begin{array}{cccccccccccccccccccccc}   
            \zeta^{'}\otimes \zeta^{'}a(t,x,\zeta)    & \zeta^{'}\otimes \zeta^{'} c(t,x,\zeta) \\[0.3cm]        
            \zeta^{'}\otimes \zeta^{'} d(t,x,\zeta)  &  \zeta^{'}\otimes \zeta^{'} b(t,x,\zeta)\\
   \end{array}  
  \right).
\end{equation}
\noindent Here $a(t,x,\zeta)$ and $b(t,x,\zeta)$ are positive measures, while
$ c(t,x,\zeta)$ and $d(t,x,\zeta)$ are complex measures such that $c=\bar d$, all supported in $[\lbrace \zeta_0 =0\rbrace \cup \lbrace \zeta ' =0\rbrace ]\cap \lbrace \zeta_1 \zeta_2\zeta_3 =0\rbrace$.\par    
\noindent They satisfy the transport system (propagation property)\par
\begin{equation}\label{ca1prooo}
\left\{
\begin{array}{ccccccccccc}
\displaystyle
 \mid \zeta '\mid^2 (-\Dis\frac{\partial a}
{\partial t}-2a )  -\Dis\sum_{l'=1}^{3}\partial_{l^{'}}.[ Tr (\zeta^{'}\otimes \zeta ' \partial^{l^{'}}{\bf E}) {c}]= 2 Re \ Tr \mu_{{uf}_{11}},\\[0.3cm] 
\Dis\sum_{l'=1}^{3}\partial_{l^{'}}.[ Tr (\zeta^{'}\otimes \zeta ' \partial^{l^{'}}{\bf E}) {a}] - \mid\zeta^{'}\mid^2  \Dis\frac{\partial c}{\partial t}= 2 Re \ Tr \mu_{{uf}_{12}},\\[0.3cm] 
 \mid \zeta '\mid^2 (-\Dis\frac{\partial b}
{\partial t}) +\Dis\sum_{l'=1}^{3}\partial_{l^{'}}.[ Tr (\zeta^{'}\otimes \zeta ' \partial^{l^{'}}{\bf E}) {d}]= 2 Re\ Tr \mu_{{uf}_{22}},\\[0.3cm]
-\Dis\sum_{l'=1}^{3}\partial_{l^{'}}.[ Tr (\zeta^{'}\otimes \zeta ' \partial^{l^{'}}{\bf E}) {b}] + \mid\zeta^{'}\mid^2  (\Dis\frac{\partial d}{\partial t}-2d)= 2Re \ Tr \mu_{{uf}_{21}}.\\[0.3cm] 
\end{array}
\right.
\end{equation}
\noindent Above, a derivative with an upper (resp. lower) index denotes a derivative wrt. variable $\zeta '$ (resp. $x$). $\bf E =E(\zeta ')$ is the constant $3\times 3$ matrix whose action is given by $E.\alpha =\zeta '\wedge \alpha$, for all $\alpha\in \R^3$. Finally, $\mu_{uf}$ is the $6\times 6$ matrix correlating the sequences $u^{\Dis\veps}$ and $f^{\Dis\veps}$, written with blocks of size $3\times 3$. In particuliar, it is zero if at least $f^{\Dis\veps}$ is strongly convergent to $0$.\par
\noindent Finally, one has the following constraint
$$\zeta^2_1 \partial_{x_1} a + \zeta^2_2 \partial_{x_2} a +\zeta^2_3 \partial_{x_3} a =2 Re \ Tr \mu_{u\tilde\varrho_{11}}$$
and similarly for $b$, $c$ and $d$, where $\mu_{u\tilde\varrho}$ is the $6\times 6$ matrix correlating the sequence $u^{\Dis\veps}$ with the sequence $\tilde\varrho^{\Dis\veps} \equiv (\varrho^{\Dis\veps},0,0,0,0,0)^t$.\par
\end{theorem}

\begin{theorem}
\label{theoreme2}{\bf Non constant coefficient-scalar case}\\ 
\noindent Assume (\ref{lois-const}), (\ref{hyp-conv-donnees}), (\ref{hyp-conv-u}) and (\ref{hyp2}). Let the dispersion matrix (see formula (\ref{definition-L}) below) be $L(x,\zeta)=\Dis\sum_{j=1}^{3}A_{0}^{-1}(x) \zeta_j A^j$, which 
has three eigenvalues, each with fixed multiplicity two, for $\zeta '\neq 0$ and given by
$$\omega_{0}=0\hspace{0.2cm}, \omega_{+}=v|\zeta^{'}|\hspace{0.2cm}, 
\omega_{-}=-v|\zeta^{'}|.$$ 
\noindent Then the matrix $P(x,\zeta)=A_{0}(\zeta_{0} {\bf Id}+ L(x,\zeta))$ has also the 
following three eigenvalues
$$\omega_{0}=\zeta_{0}\hspace{0.2cm}, \omega_{+}=\zeta_{0}+v|\zeta^{'}|\hspace{0.2cm}, 
\omega_{-}=\zeta_{0}-v|\zeta^{'}|,$$ 
\noindent each with fixed multiplicity two, where $v(x)=\Dis\fa{1}{\sqrt{\epsilon(x)\eta(x)}}$ 
is the propagation speed.\par
\noindent Using the propagation basis and the eigenvector basis introduced in (\ref{vecteurs-propres1}) and (\ref{vecteurs-propres2}), it follows that the H-measure $\scorpio$ associated to a suitable subsequence of $u^{\Dis\veps}$ can be expressed as:
$$
\scorpio = \left(\begin{array}{cccccccccc} \scorpio_{11} & \scorpio_{12} \\ \scorpio_{21} &\scorpio_{22}\\ \end{array}\right)
$$
where $\scorpio_{ij}$ are $3\times 3$ matrix valued measures. Furthermore, one has
$$
\left\{
\begin{array}{cccccc}
\displaystyle
\scorpio_{11}=
       \fa{1}{\epsilon} [(\hat{\zeta '}\otimes \hat{\zeta '})a_{0}+
\fa{1}{2}(z^{1}\otimes z^{1}) a_{+}+\fa{1}{2}(z^{2}\otimes z^{2}) b_{+}
+\fa{1}{2}(z^{1}\otimes z^{1})a_{-}+\fa{1}{2}(z^{2}\otimes z^{2}) b_{-}]\\[0.3cm]
\scorpio_{12}=\fa{v}{2}
[(z^{1}\otimes z^{2})a_{+}-
(z^{2}\otimes z^{1}) b_{+}-(z^{1}\otimes z^{2})a_{-}+
(z^{2}\otimes z^{1}) b_{-}]\\[0.3cm]
\scorpio_{21}=
\fa{v}{2}[(z^{2}\otimes z^{1}) a_{+}-(z^{1}\otimes z^{2})b_{+}-
(z^{2}\otimes z^{1})a_{-}+(z^{1}\otimes z^{2})b_{+}]\\[0.3cm]
\scorpio_{22}= 
\fa{1}{\mu}[(\hat{\zeta '}\otimes \hat{\zeta '})b_{0}+\fa{1}{2}(z^{2}\otimes z^{2}) a_{+}+\fa{1}{2}(z^{1}\otimes z^{1})b_{+}+
\fa{1}{2}(z^{2}\otimes z^{2})a_{-}+\fa{1}{2}(z^{1}\otimes z^{1})b_{-} ]\;.
\end{array}
\right.\;
$$
using notations given by (\ref{vecteurs-propres1}). Above $a_0$, $b_0$, $a_\pm$ and $b_\pm$ are all scalar positive measures supported in the set $[\{\zeta_0 = 0\}\cup \{\zeta_0 =\pm v\mid\zeta '\mid \}]\cap \{\zeta_1 \zeta_2\zeta_3 =0\}$. Finally, one has the following propagation type system
$$
\left\{\matrix{\displaystyle 
-\Dis\veps (x) \partial_t \scorpio_{11} +\zeta_0\sum^3_{l'=1} \partial_{l'}\veps (x) \partial^{l'}\scorpio_{11} - 2\sigma \scorpio_{11} - \sum^{3}_{l'=1} \partial^{l'} {\bf E} .\partial_{l'} \scorpio_{12} = 2 Re \mu_{{uf}_{11}}\cr
\displaystyle -\eta (x) \scorpio_{12} +\zeta_0 \sum^3_{l'=1} \partial_{l'} \eta (x) \partial^{l'} \scorpio_{12} +\sum^3_{l' =1} \partial^{l'} {\bf E} \partial_{l'}\scorpio_{11} = 2 Re \mu_{{uf}_{12}} \cr
\displaystyle 
-\Dis\veps (x) \partial_t \scorpio_{21} +\zeta_0\sum^3_{l'=1} \partial_{l'}\veps (x) \partial^{l'}\scorpio_{11} - 2\sigma \scorpio_{21} - \sum^{3}_{l'=1} \partial^{l'} {\bf E} .\partial_{l'} \scorpio_{22} = 2 Re \mu_{{uf}_{21}}\cr
\displaystyle -\eta (x) \scorpio_{22} +\zeta_0 \sum^3_{l'=1} \partial_{l'} \eta (x) \partial^{l'} \scorpio_{22} +\sum^3_{l' =1} \partial^{l'} {\bf E} \partial_{l'}\scorpio_{21} = 2 Re \mu_{{uf}_{22}} \\
}\right.
$$
where we are using same notations as in Theorem \ref{theoreme1} for the right hand side.
\end{theorem}

\noindent \section{Some basic facts on H-measures}
\noindent In this Section, we recall some results from the H-measures theory, taking the presentation of Tartar \cite{tar}. However, this is also similar to the exposition of G\'erard \cite{Gér1, Gér3}, relying upon 
Hormander \cite{Horm1},  \cite{Horm2}.
\noindent \begin{Definition}\label{dhme} Let $\Omega$ be an open set of $\R^{\;n}$ and let $ u^{\displaystyle\veps}$ be 
a sequence of functions defined in $\R^{\;n}$ with values in $\R^{\;p}$. 
We assume that $ u^{\displaystyle\veps}$ converges weakly to zero in $ (L^2(\R^{n}))^p $. 
Then after extracting a subsequence (still denoted 
by $\Dis\veps $), there exists a family of complex-valued Radon measures 
$(\scorpio_{\;ij}(x,\zeta ))_{1\leq i,j \leq p }$ on $\R^{\;n}\times S^{\;n-1}$,  
such that for every functions $\phi_1,\phi_2$ in $C_{0}(\R^{\;n})$, 
the space of continuous functions converging to zero at infinity, and for 
every function $\psi $ in $C(S^{\;n-1})$, the space of continuous functions on 
the unit sphere $S^{\;n-1}$ in $\R^{n}$, one has
\begin{equation} \label{defini 1}
\left\{
\begin{array}{ccccccccc}
<\scorpio_{\;ij},\phi_1\bar{\phi_2}\otimes \psi >=\Dis\int_{\R^{n}}
\Dis\int_{S^{\;n-1}}\phi_1\;\bar{\phi_2}
\psi(\zeta/|\zeta |)\scorpio_{\;ij}(x,\zeta ) dx d\zeta \\[0.3cm]
=\Dis\lim_{\Dis\veps\longrightarrow 0}
\int_{\R^{n}} [F(\phi_1u_{i}^{\displaystyle\veps})(\zeta)]
\overline{[F(\phi_2u_{j}^{\displaystyle\veps})(\zeta)]}
\psi(\zeta/|\zeta |) d\zeta.
\end{array}
\right.
\end{equation}
\end{Definition}    
               
\noindent Above, $ F$ denotes the Fourier transform operator defined in 
$L^2(\R^{n})$, for an integrable function $ f$ as 
$[F(f)](\zeta)=\Dis\int_{\R^{n}}f(x)\; e^{-2\pi i x.\zeta}  dx $, 
while $\bar{F}$ is the inverse Fourier transform defined as 
$\overline{[F(f)]}(x)=\Dis\int_{\R^{n}}f(\zeta )\; e^{2\pi i x.\zeta}  d\zeta $. $\bar{z}$ denotes the complex conjugate of the complex number $ z$.\par
\noindent The matrix valued measure $\scorpio=(\scorpio_{\;ij})_{1\leq i \leq p }$ is 
called the H-measure associated with the extracted subsequence 
$ u^{\displaystyle\veps}$ .\par
\noindent\begin{Remark} 
\noindent H-measures are hermitian and 
non-negative matrices in the following sense  

\begin{equation}\label{2.20}
\left\{
\begin{array}{ccc}
\scorpio_{ij}=\overline{\scorpio_{ji}} \hspace{1cm}\mbox{and} \\[0.3cm]
\Dis\sum_{i,j=1}^{p}\scorpio_{ij} \phi_i\overline{\phi_j}\geq 0\;\; 
\mbox{for all}\;\; \phi_1\;\phi_2\;...........\phi_n\in C_{0}(\R^n),  
\end{array}
\right.
\end{equation}

\noindent and it is clear that the H-measure for a strongly convergent sequence is zero.\par   
\end{Remark}
\noindent Although we consider the scalar case for all properties of H-measures, all the following facts 
are easily extended to the vectorial case.\par
\noindent\begin{Definition}
\noindent Let $a\in C(S^{\;n-1})$, $ b \in C_{0}(\R^{n})$. We associate with $ a $ the linear continuous 
operator $ A $ on $L^2(\R^{n})$ defined by

\begin{equation}\label{2.3}
F(Au)(\zeta )=a(\zeta /|\zeta |) F(u)(\zeta )\hspace{0.2cm} a.e.
\hspace{0.2cm} \zeta \in \R^{n} 
\end{equation}

\noindent and with $ b $ we associate the operator

\begin{equation} \label{2.4}
Bu(x)=b(x)u(x)\hspace{0.2cm} a.e.\hspace{0.2cm}  x\in \R^{n}\;.
\end{equation}

\noindent A continuous function $ P $ on $\R^{\;n}\times S^{\;n-1}$ with values 
in $\R$ is called an admissible symbol if it  can be written as

\begin{equation}\label{2.5}
P(x,\zeta )=\Dis\sum_{n=1}^{+\infty} b_{n}(x)\otimes 
a_{n}(\zeta )=\Dis\sum_{n=1}^{+\infty} b_{n}(x)a_{n}(\zeta)  
\end{equation}

\noindent where $a_{n}$ are continuous functions on $S^{\;n-1} $ 
and $ b_{n}$ are continuous bounded functions converging to zero at 
infinity on $\R^{\;n} $ with

\begin{equation} \label{2.6}
\Dis\sum_{n=1}^{+\infty} \max_{\Dis \zeta } |a_{n}(\zeta)| \max_{\Dis x } 
|b_{n}(x)|<\infty \;.
\end{equation} 
\end{Definition}

\noindent An operator $L$ with symbol $P$ is defined by:\par 
\noindent 1) $L$ is linear continuous on $L^2(\R^{n})\;.$\par
\noindent 2) $P$ is an admissible symbol with a decomposition (\ref{2.5}), satisfying 
(\ref{2.6}).\par    
\noindent 3) $L$ can written as the following form
$$L=\Dis\sum_{n=1}^{n} A_{n} B_{n}+\mbox{compact operator}$$

\noindent where $ A_{n}\;, B_{n} $ are the operators associated with 
$ a_{n}\;,b_{n}$ as in (\ref{2.3}) and (\ref{2.4}).\par   

\noindent With notations as in (\ref{2.3}) and (\ref{2.4}), one can 
show that the operator 
$ C:=AB-BA $ is a compact operator from $L^2(\R^{n})$ into itself. 
\noindent Denote by $X^{m}(\R^{n})$ the space of functions $ v $ 
with derivatives up to order $ m $ belonging to the image by the Fourier 
transform of the space $ L^1(\R^{n})$ i.e. $(F(L^1(\R^{n})))$, equipped 
with the norm
$$||v|||_{X^{m}}=\int_{\R^{n}}(1+|2\pi \zeta |^m)|F(v)(\zeta )| \;d\zeta \;.$$  

\noindent Then, if  $ A $ and $ B$ are operators with symbols 
$ a $ and $ b$ as in (\ref{2.3}) and (\ref{2.4}), satisfying one the
following conditions \par
\noindent 1)  $ a\in C^1(S^{\;n-1})$ and $ b\in X^1(\R^{n})\;,$ \par
\noindent 2) $ a \in X_{loc}^{1}(\R^{n} \setminus \{0 \})$ and $ b\in 
C_{0}^{1}(\R^{n}) $,\par 
\noindent it follows that the operator $ C=AB-BA $ is a 
continuous operator from $ L^2(\R^{\;n})$ into $ H^{1}(\R^{\;n}) $ and 
extending $ a $ to be homogeneous of degree zero on $ \R^{\;n} $,  
then $ \nabla C=\Dis\frac{\partial }{\partial x_i} (AB-BA) $ has the symbol

\begin{equation} \label{2.7}
(\nabla_{\Dis\zeta } a. \nabla_{\Dis x} b )\zeta =\zeta_{i}
\sum_{ k=1}^{n} \frac{\partial a}{\partial \zeta_{k}}\;
\frac {\partial b}{\partial x_{k}}\;.
\end{equation}    
\smallskip
\noindent The main results of H-measures theory are given by the 
next two results\par 

\begin{theorem}\label{loca}{\bf Localisation property}
Let $u^{\;\displaystyle\veps}$ be a sequence converging weakly to zero in 
$ (L^2(\R^{n}))^p$ and let $\scorpio$ be the H-measure associated 
to $u^{\displaystyle\veps}$. Assume that one has the balance relation   \\
$$\Dis\sum_{k=1}^{n} \frac{\partial }
{\partial x_{k}}(A^{k} u^{\displaystyle\veps})\longrightarrow  0  \hspace{0.2cm}
(H_{loc}^{-1}(\Omega ))^p \hspace{0.2cm}\mbox{strongly}, $$  
\noindent where $ A^{k}$ are continuous matrix valued functions on 
$\Omega \subset \R^{\;n}$. Then, on $\Omega \times S^{\;n-1}$, one has 

\begin{equation}\label{2.10}
P(x,\zeta)\scorpio \equiv (\Dis\sum_{k=1}^{n}\zeta_k A^k(x))\scorpio =0 \;.         
\end{equation}       
\end{theorem}   

\noindent This result shows that the support of the $H$-measure $\scorpio $ is 
contained in the (characteristic) set  $$\{(x,\zeta) \in \Omega \times S^{\;n-1}\;,\;\; 
\mbox{det}\;  P(x,\zeta)=0\}\;.$$

\noindent \begin{theorem}\label{propaga}
{\bf Propagation property for symmetric systems}\label{propagation-theorie}
\noindent Let be given matrix valued functions $A^{k}$ in the class 
$ C_{0}^1(\Omega )$. Assume that the pair of sequences $(u^{\displaystyle\veps},f^{\displaystyle\veps})$ 
satisfies the symmetric system 
\begin{equation} \label{2.8}
\Dis\sum_{k=1}^{n} A^k \frac{\partial u^{\displaystyle\veps}}
{\partial x_{k}}+ B u^{\displaystyle\veps}= 
f^{\displaystyle\veps}
\end{equation}   

\noindent and that both sequences $(u^{\displaystyle\veps})$,  
$(f^{\displaystyle\veps})$ converge weakly to zero in 
$L^2(\Omega )^p$. Then the H-measure $\mu $ associated to the sequence $(u^{\;\displaystyle\veps},
f^{\;\displaystyle\veps})$ and given under the form

\begin{equation}\label{forme-bloc-theorie}
\mu =\left(
           \begin{array}{cccccccccccccccccccccc}   
            {\bf\mu_{11}}  & {\bf\mu_{12}}\\[0.3cm]        
            {\bf\mu_{21}}  &  {\bf\mu_{22}} \\
   \end{array}  
  \right)
\end{equation}

\noindent satisfies the equation 

\begin{equation}\label{2.11}
<{\bf\mu_{11}},\{P,\psi\}+\psi \sum^n_{k=1}\partial_{k}A^k-2\psi S>=<2\Re( 
\mbox{Tr}{\bf\mu_{12}}),\psi>
\end{equation} 

\noindent for all smooth functions $\psi (x,\zeta )$. Here $ S:=1/2(B+B^{\;*})$ is the hermitian part of the matrix $B$ and $\{P,\psi\}$ is the Poisson bracket of $ P $ and 
$\psi $, i.e.\par

\begin{equation}\label{2.12}
\{P,\psi\}=\partial^{\;l} P \partial_{\;l}\psi -\partial^{\;l}\psi 
\partial_{\;l} P \equiv 
\Dis\sum_{l=1}^{n} (\frac{\partial P}{\partial \zeta_l} \frac{\partial 
\psi }{\partial x_l}-
\frac{\partial \psi }{\partial \zeta_l}\frac{\partial P}{\partial x_l}).
\end{equation} 
\end{theorem}

\section{ Applications to Maxwell's system}  
\noindent This section is devoted to the proofs of our main results stated in the Introduction.\par
\subsection{ Proof of Theorem \ref{theoreme1}: Constant coefficient case}

\noindent This case corresponds to the assumption (\ref{hyp1}), that is all the matrices $\ddot{\epsilon} \;,\ddot{\eta} \;,\mbox{and}\; 
\ddot{\sigma}$ are the identity matrix, i.e\par
\begin{equation}
\ddot{\epsilon} =\ddot{\sigma} =\ddot{\eta} =\left(
                  \begin{array}{cccccccccccccccccccccc}   
                    1 &  0  & 0 \\
                     0 &  1  & 0  \\
                     0 &  0  & 1   \\
   \end{array}  
  \right)= (Id)_{3\times 3}\;.
\end{equation}

\noindent In this case, system (\ref{maxwel}) can be rewritten as       

\begin{equation}\label{maxwell-hyp1}
\left\{
\begin{array}{cccccl} 
i)&\Dis\frac{\partial D^{\Dis\veps}}{\partial t}(x,t) +
E^{\Dis\veps}(x,t)& =&\mbox{rot} H^{\Dis\veps}(x,t)+ F^{\Dis\veps}(x,t),&\\[0.3cm]  
ii)&\displaystyle\frac{\partial H^{\Dis\veps}}{\partial t}(x,t)& 
          =&-\mbox{rot}E^{\Dis\veps}(x,t)+ G^{\Dis\veps}(x,t),&\\[0.3cm]  
iii)&\Dis\mbox{div} H^{\Dis\veps}(x,t)&=&0,& \\[0.3cm]            
iv)&\Dis\mbox{div} E^{\Dis\veps}(x,t)&=&\rho^{\Dis\veps}(x,t)\;.&\\[0.3cm]     
\end{array}
\right.
\end{equation}  
\noindent Recalling the notation of the Introduction, it follows that Maxwell's system (\ref{maxwell-hyp1}) can be written as
\begin{equation}\label{systeme1-constant}
\Dis\sum_{i=0}^{3} A^i \frac{\partial u^{\displaystyle\veps}}{\partial x_{i}}+
C u^{\displaystyle\veps}= 
f^{\displaystyle\veps}
\end{equation}    
and
\begin{equation}\label{systeme2-constant}
\Dis\sum_{i=1}^{3} B^i \frac{\partial u^{\displaystyle\veps}}{\partial x_{i}}= \tilde\varrho^{\Dis\veps}.
\end{equation}
\noindent Here

\begin{equation}\label{A0-constant}
 A^{\;0} =\left(
                  \begin{array}{cccccccccccccccccccccc}   
                   \ddot{\epsilon}   &  \bf{0}  \\
                     \bf{0}      &  \ddot{\eta}     \\

  \end{array}  
  \right)=\left(
                  \begin{array}{cccccccccccccccccccccc}   
                     \bf{Id}  &  \bf{0}  \\
                      \bf{0}      &      \bf{Id}     \\
\end{array}  
  \right)
\end{equation}  

\noindent and

\begin{equation}\label{elke}
  A^{\;1} =\left(
                  \begin{array}{cccccccccccccccccccccc}   
                     \bf{0} &  \bf{Q_1}^{t}  \\
                     \bf{Q_1}    &  \bf{0}       \\
  \end{array}  
  \right)\;,
 A^{\;2} =\left(
                  \begin{array}{cccccccccccccccccccccc}   
                     \bf{0} &  \bf{Q_2}^{t}  \\
                     \bf{Q_2}    &  \bf{0}     \\
  \end{array}  
  \right)\;,
 A^{\;3} =\left(
                 \begin{array}{cccccccccccccccccccccc}   
                     \bf{0} &  \bf{Q_3}^{t}  \\
                     \bf{Q_3}    &  \bf{0}   \\
  \end{array}  
  \right)\;.
\end{equation}

\noindent The constant antisymmetric matrices ${\bf Q_k} \;, 1\leq k\leq 3$ 
are given by

\begin{equation}
 \bf{Q_1} =\left(
                  \begin{array}{cccccccccccccccccccccc}   
                     0 &  0  & 0 \\
                     0 &  0  & -1  \\
                     0 &  1  & 0   \\
   \end{array}  
  \right)\;,
\bf{Q_2} =\left(
                  \begin{array}{cccccccccccccccccccccc}   
                     0 &  0  & 1 \\
                     0 &  0  & 0  \\
                    -1 &  0  & 0   \\
   \end{array}  
  \right)\;,
\bf{Q_3} =\left(
                  \begin{array}{cccccccccccccccccccccc}   
                     0 &  -1  &  0 \\
                     1 &   0  &  0  \\
                     0 &   0  &  0   \\
   \end{array}  
  \right)
\end{equation}

\noindent the matrix $ C$ and $f^{\Dis\veps}$ by

\begin{equation}\label{C-constant}
C= \left(
   \begin{array}{cccccccccccccccccccccc}   
    \ddot{\sigma} & \bf{0}\\
     \bf{0}     & \bf{0} \\
   \end{array}  
  \right)=
\left(
   \begin{array}{cccccccccccccccccccccc}   
    \bf{Id} & \bf{0}\\
     \bf{0}     & \bf{0} \\
   \end{array}  
  \right)\;, 
f^{\Dis\veps}=\left(
   \begin{array}{cccccccccccccccccccccc}   
    F^{\Dis\veps} \\
     G^{\Dis\veps}\\
   \end{array}  
  \right)\;.
\end{equation}
Matrices $B^i$, $i=1,2,3$, are given by
$$B^i = \left(
   \begin{array}{cccccccccccccccccccccc}   
    \beta^i & \bf{0}\\
     \bf{0}     & \beta^i \\
   \end{array}  
  \right)$$
where the $3\times 3$ matrices $\beta^i$ are such that
$$\beta^i_{kl} =0 \mbox{ except for } \beta^i_{ii} =1$$
Finally we have denoted $\tilde \varrho^{\Dis\veps} \equiv (\varrho^{\Dis\veps},0,0,0,0,0)^t$.\par
\noindent Denote the H-measure corresponding to (a subsequence of) the 
sequence $u^{\Dis\veps}$ by

\begin{equation}\label{3.9}
\scorpio=\left(
           \begin{array}{cccccccccccccccccccccc}   
            \bf{\nu_{e}}   &  \bf{\nu_{em}}\\[0.3cm]        
            \bf{\nu_{me}}  &  \bf{\nu_{m}}  \\
   \end{array}  
  \right)\;.
\end{equation}

\noindent The measure $\scorpio$ is a $2\times 2$ block matrix measure, each block being of size $3\times 3$.\par
\noindent In the following, let $x_0 =t$, $\tilde x =(x_0, x)$, $x=(x_1,x_2,x_3)$. We let $\zeta$ denote the dual variable to $\tilde x$, with $\zeta =(\zeta_0 ,\zeta ')$, $\zeta ' =(\zeta_1, \zeta_2,\zeta_3)$.\par
\noindent To state the localisation property (\ref{loca}), we need first to express the symbol of the differential operator appeating in (\ref{systeme1-constant}), for which one has\par

$$ P(x,\zeta) \equiv \Dis\sum_{j=0}^{3}\zeta_j A^j(x) \;$$

\noindent  and thus

$$
P(x,\zeta)=\zeta_0 A^0(x)+\zeta_1 A^1(x) +\zeta_2 A^2(x) +\zeta_3 A^3(x)=
$$

$$=
   \zeta_0 \left(
               \begin{array}{cccccccccccccccccccccc}   
                  \bf{Id}    & \bf{0} \\
                     \bf{0}  &  \bf{Id} \\
  \end{array}  
  \right)+
 \zeta_1 \left(
               \begin{array}{cccccccccccccccccccccc}   
                  \bf{0} & \bf{Q_1^t} \\
                  \bf{Q_1}   &  \bf{0} \\
  \end{array}  
  \right)
+
\zeta_2 \left(
               \begin{array}{cccccccccccccccccccccc}   
                  \bf{0} & \bf{Q_2^t} \\
                  \bf{Q_2}   &  \bf{0} \\
  \end{array}  
  \right)$$
$$+
\zeta_3 \left(
               \begin{array}{cccccccccccccccccccccc}   
                  \bf{0} & \bf{Q_3^t} \\
                   \bf{Q_3}   &  \bf{0} \\
  \end{array}  
  \right)=
        \left(
               \begin{array}{cccccccccccccccccccccc}   
                  \zeta_0\bf{Id} & -\bf{E} \\
                  \bf{E}          &  \zeta_0\bf{Id}  \\
  \end{array}  
  \right)
$$

\noindent where

\begin{equation}\label{matr}
\bf{E} \equiv \left(
               \begin{array}{cccccccccccccccccccccc}   
                   0       & -\zeta_3   & \zeta_2   \\
                  \zeta_3  & 0          & -\zeta_1\\
                  -\zeta_2 & \zeta_1   & 0  \\
  \end{array}  
\right) \;.
\end{equation}
\noindent Clearly $ \bf{E}$ is antisymmetric (i.e. $ \bf{E^{t}}=-\bf{E} $), 
so that $P$ is a symmetric matrix.\par
\noindent Using the localisation property, it follows

\begin{equation}\label{localisation-constant}
P\scorpio=\left(
               \begin{array}{cccccccccccccccccccccc}   
                  \zeta_0\bf{Id} & -\bf{E}\\[0.3cm]        
                   \bf{E}      &  \zeta_0\bf{Id}  \\
                \end{array} 
 \right)
           \left(
           \begin{array}{cccccccccccccccccccccc}   
            \bf{\nu_{e}}   &  \bf{\nu_{em}}\\[0.3cm]        
            \bf{\nu_{me}}  &  \bf{\nu_{m}}  \\
   \end{array}  
  \right)=\bf{0}
\end{equation}

\noindent and thus\par 

\begin{equation}\label{3.11}
\left\{
\begin{array}{cccccccc}
\displaystyle 
&1)&&\zeta_0\bf{Id}.\bf{\nu_{e}}+\bf{E^t}.\bf{\nu_{me}} =&&0,&\\ 
&2)&&\zeta_0\bf{Id}.\bf{\nu_{em}}+\bf{E^t}.\bf{\nu_{m}}  =&&0,&\\  
&3)&&\bf{E}.\bf{\nu_{e}}+\zeta_0\bf{Id}.\bf{\nu_{em}} =&&0,&\\
&4)&&\bf{E}.\bf{\nu_{em}}+\zeta_0\bf{Id}.\bf{\nu_{m}}   =&&0.&\\
\end{array}
\right.
\end{equation}
First note that from (\ref{localisation-constant}), since (see also next subsection) $P$ has $\zeta_0$, $\pm \mid \zeta '\mid$ as eigenvalues, that $\scorpio$ is supported in $\lbrace \zeta_0 =0\rbrace \cup \lbrace \zeta ' =0 \rbrace$. Then, one has
\noindent\begin{Lemma}\label{vecttt}
\noindent The H-measure $\scorpio$ can be written under the form
\begin{equation}\label{decomposition-constant}
\scorpio=\left(
               \begin{array}{cccccccccccccccccccccc}   
                \zeta^{'}\otimes \zeta^{'}a(t,x,\zeta)  & 
		\zeta^{'}\otimes \zeta^{'} {c}(t,x,\zeta)  \\[0.3cm]        
                \zeta^{'}\otimes \zeta^{'}d(t,x,\zeta)      &  
	        \zeta^{'}\otimes \zeta^{'}b(t,x,\zeta)    \\
                \end{array} 
 \right)\;
\end{equation}
\noindent where $a$, $b$ are scalar positive measures, $c$ and $d$ are scalar complex measures such that $\bar c =d$, all supported in $\lbrace \zeta_0 =0 \rbrace \cup \lbrace \zeta ' =0 \rbrace$. 
\end{Lemma}
\smallskip 
\noindent {\bf Proof of Lemma \ref{vecttt}}\par
\smallskip
\noindent Multiplying (\ref{3.11}-1) par $ \zeta_0 $, (\ref{3.11}-3) par $\bf{E^{t}}$ 
and substracting the results, one has 

\begin{equation}\label{3.12}
(\zeta_0 ^2\bf{Id} + \bf{E^2}). \bf{\nu_{e}}=0\;.
\end{equation}

\noindent We  discuss the following distinct cases\par
\smallskip 
\smallskip 
\noindent {\bf i)} case $\zeta_0=0$. Then note that $\zeta^{'}\neq 0 \;$ since $\zeta$ belongs to the unit sphere of $\R^4$. From (\ref{3.11}), we have ${\bf E}.\bf{\nu_{e}}=0 $. \par

\noindent Then, we use the following Lemma 

\noindent \begin{Lemma}\label{mavfor}

\noindent If $ {\bf E}.A=0$, then the matrix 
$ A $ has the form $A= \zeta^{\;'}\otimes a $, for some vector 
$a\in\R^{\;3}$.\par 
\end{Lemma}

\noindent {\bf Proof of Lemma \ref{mavfor}}\par

\noindent We denote the columns of the matrix $A$ by the vectors 
\begin{equation} \label{eqmatz}
A=[ \vec{a}_{1} \;\;\vec{a}_{1} \;\;\vec{a}_{1}]\;.
\end{equation}

\noindent But as ${\bf E}.A=0$, we get that 
\begin{equation} 
[E \vec{a}_{1} \;\;E \vec{a}_{1} \;\;E \vec{a}_{1}]=0\;
\end{equation}
\noindent or 
\begin{equation}\label{indmat} 
E \vec{a}_{i} =0\;\;\;, i=1,2,3\;.
\end{equation}

\noindent For $i=1$, and similarly for the other cases, 
\smallskip 
\smallskip 
\noindent $a_{1}= \left(
                \begin{array}{cccccccccccccccccccccc}   
                     v_{1}\\
                     v_{2}\\
                     v_{3}
                 \end{array}
\right )\in \R^{3}$. Then from (\ref{matr}), (\ref{indmat}), we get that

\begin{equation}\label{matrtghlu}
\left(
               \begin{array}{cccccccccccccccccccccc}   
                   0       & -\zeta_3   & \zeta_2   \\
                  \zeta_3  & 0          & -\zeta_1\\
                  -\zeta_2 & \zeta_1   & 0  \\
  \end{array}  
\right) 
\left(
                \begin{array}{cccccccccccccccccccccc}   
                     v_{1}\\
                     v_{2}\\
                     v_{3}
                 \end{array}
\right )=\left(
               \begin{array}{cccccccccccccccccccccc}   
                   -v_{2}\zeta_3 + v_{3}\zeta_2   \\
                    v_{1}\zeta_3 - v_{3}\zeta_1\\
                   -v_{1}\zeta_2 +v_{2} \zeta_1   \\
  \end{array}  
\right)=\zeta^{'}\otimes \vec{a}_{1}=0\;
\end{equation}

\noindent which implies that $\zeta^{'}// \vec{a}_{1}$ and thus $\vec{a}_{1}=c_{1}\zeta^{'}$ , for some 
constant $c_{1}\in\R$.

\noindent Thus all in all, and for $i=1,2,3$, all the columns of the matrix $A$ are 
parallel to the vector $\zeta^{'}=(\zeta_1,\zeta_2,\zeta_3)$, so we can write that $\vec{a}_{i}=c_{i}\zeta^{'}$, and 
by arranging these numbers $c_{i}$ as components of the vector $a$, we get that   
\begin{equation} 
A=\vec{a}\otimes \zeta^{'}.
\end{equation}
\smallskip 
\smallskip 
\noindent {\bf End of the proof of Lemma \ref{vecttt}}\par
\noindent Using Lemma \ref{mavfor}, we can conclude that
${\bf \nu_{e}}=\zeta^{'}\otimes \zeta^{'} a(t,x,\zeta^{'})$ and thus  
the blocks of the matrix $H$-measure, which satisfy system (\ref{3.11}), are such that

\begin{equation}\label{2.18}
\left\{
\begin{array}{ccccccc}
\displaystyle
&{\bf\nu_{e}}= \zeta^{'}\otimes \zeta ^{'} a(t,x,\zeta),&\\ [0.3cm] 
&{\bf\nu_{m}} = \zeta{'}\otimes \zeta^{'}b(t,x,\zeta),&\\[0.3cm]  
&{\bf\nu_{em}} =\zeta^{'}\otimes \zeta^{'}c(t,x,\zeta),&\\[0.3cm]  
&{\bf\nu_{me}} =\zeta^{'}\otimes \zeta^{'}d(t,x,\zeta )=
{\bf\bar{\nu}_{em}}= \zeta^{'}\otimes \zeta^{'}\bar{c}(t,x,\zeta)&\\ 
\end{array}
\right.
\end{equation}
\smallskip 

\noindent where $a(t,x,\zeta)$ and $b(t,x,\zeta)$ are real 
positive measures, while $d(t,x,\zeta)$ and $c(t,x,\zeta)$ are scalar complex 
measures, such that $c=\bar{d}$.\par
\smallskip
\smallskip
\noindent {\bf ii)} case $\zeta_0\neq 0\;\mbox{but}\; \zeta^{\;'}=
(\zeta_1,\zeta_2,\zeta_3)=0$. \par

\noindent From (\ref{3.12}) one has $\zeta_0^{\;2}\bf{Id}.\bf{\nu_{e}}=0$ 
(recall that $\zeta^{\;'}=0 $ implies that $ \bf{E}=0$), and in this case we have 
$\zeta_0^{\;2}.\bf{\nu_{e}}=0 $. As $\zeta_0\neq 0$, we get 
$ \bf{\nu_{e}}=0$. Repeating the same steps, for the other equations of the system (\ref{3.11}), we 
get that $ {\bf\nu_{e}}={\bf\nu_{m}}={\bf\nu_{em}}=\overline{{\bf\nu_{me}}}=0$ and finally, one has

\begin{equation}\label{3.15}
\scorpio=\left(
               \begin{array}{cccccccccccccccccccccc}   
                  \bf{0} &  \bf{0} \\[0.3cm]        
                  \bf{0} &  \bf{0}  \\
                \end{array} 
 \right).
\end{equation}
\smallskip
\smallskip

\noindent {\bf iii)} case $\zeta_0\neq 0\;\mbox{and}\; \zeta^{'}=
(\zeta_1,\zeta_2,\zeta_3)\neq 	0$. Using (\ref{3.12}), one has

\begin{equation}\label{3.16} 
(\zeta_0^{\;2}{\bf{Id}} + \bf{E}.\bf{E}). \bf{\nu_{e}}=0\Longleftrightarrow 
\zeta_0^{\;2}.\bf{\nu_{e}}+
\bf{E}.(\bf{E}.\bf{\nu_{e}})=0.
\end{equation}
\noindent But as $(\bf{E}.\bf{\nu_{e}})=
\vec{\zeta^{'}} \otimes {\bf{\nu_{e}}}$, one has that $ \bf{E}.
( \vec{\zeta^{'}} \otimes {\bf{\nu_{e}}})=\vec{\zeta^{'}} \otimes (\vec{\zeta^{'}} 
\otimes {\bf{\nu_{e}}})$ and from (\ref{3.16}), one gets

\begin{equation}\label{3.177}
\zeta_0^{\;2}.\bf{\nu_{e}}+
\vec{\zeta^{'}} \otimes (\vec{\zeta^{'}} \otimes {\bf{\nu_{e}}})=0.
\end{equation}


\noindent Now we shall show that $ {\bf\nu_{e}}=0 $.\par
\smallskip
\smallskip

\noindent By contradiction, we assume that $ \vec{\nu}_{\bf{e}}\neq 0 $. 
Using $(\ref{3.177})$, one has $ (1-\zeta^{'2}).{\bf\nu_{e}}=
-\vec{\zeta^{'}} \otimes (\vec{\zeta^{'}} \otimes {\bf{\nu_{e}}}) $, thus either $\zeta^{'2}=1\;\mbox{or}\; 
\zeta^{'2}< 1\;\mbox{or}\; \zeta^{'2}> 1 $.\par 
\noindent If $\zeta^{'2}=1 $, one has $ \zeta_0 = 0 $, which is a 
contradiction with $ \zeta_0 \neq 0 $. If $ \zeta^{'2}> 1 $, then one 
have $ \zeta_0^{\;2}=1-\zeta^{'2}< 0 $ which is not possible, because 
$\zeta_0^{\;2}> 0 $, and if $\zeta^{'2}< 1$ then we get again 
$\zeta_0^{\;2}=\zeta^{'2}-1<0$, which is a contradiction.\par

\noindent Thus all in all, we conclude that $ {\bf{\nu_{e}}}=0 $ and similarly for 
$ \nu_{m}={\bf\nu_{em}}=\overline{{\bf\nu_{me}}}=0$, and finally  we have also, in this case 

$$\scorpio=\left(
               \begin{array}{cccccccccccccccccccccc}   
                  \bf{0} &  \bf{0} \\[0.3cm]        
                  \bf{0} &  \bf{0}  \\
                \end{array} 
 \right).$$
ending the proof of Lemma \ref{localisation-constant}.\par
\smallskip
\noindent Let us now turn to the localisation property associated with equation (\ref{systeme2-constant}). In this case, it follows that $\scorpio$ satisfies
$$(\sum^3_{j=1} B^j\zeta_j ) \scorpio = 0.$$
\noindent Setting $B(\zeta ) =\Dis\sum_{\Dis j=1}^{3} B^{j}\zeta_{j} $, since $det\ B (\zeta ) =\zeta_1\zeta_2\zeta_3$, it follows that $\scorpio$ is supported in the set $\lbrace \zeta_1\zeta_2\zeta_3 =0\rbrace$.\par
\noindent All in all, the scalar measures $a$, $b$, $c$ and $d$ are all (also) supported in the set $\lbrace \zeta_1\zeta_2\zeta_3 =0\rbrace$.\par
\noindent Now, we wish to write down the propagation property, and for this purpose, we need to compute the 
Poisson bracket, associated with equations (\ref{systeme1-constant}) and (\ref{systeme2-constant}).\par
\noindent Letting $\psi=\psi(\tilde x,\zeta)$ be an arbitrary smooth function, recall first that the 
Poisson bracket is given by

\begin{equation}\label{3.19}
\left\{
\begin{array}{ccc}
\displaystyle
\{P,\psi\}=\sum^3_{l=0}\partial^l P \partial_l\psi -\partial^l\psi \partial_l P 
\\[0.4cm]=
\Dis\sum_{l=0}^{3} (\frac{\partial P}{\partial {\bf\zeta_l}} \frac{\partial \psi }
{\partial {\tilde x_l}}-
\frac{\partial \psi }{\partial {\bf \zeta_l}}\frac{\partial P}{\partial {\tilde x_l}}).
\end{array}
\right.
\end{equation}
\noindent Recall that a derivative with an upper (resp. lower) index denotes a derivative wrt. variable $\zeta$ (resp. $\tilde x$).\par
\noindent In our case, starting with (\ref{systeme1-constant}), we have

\begin{equation}\label{3.20}
\{P,\psi\}= \left(
             \begin{array}{cccccccccccccccccccccc}   
                {\bf Id} {\partial_t \psi}  &   
		-(\Dis\sum_{l'=1}^{3} \partial^{l^{'}} \bf{E}.\partial_{l^{'}}\psi )    \\[0.3cm]        
               (\Dis\sum_{l'=1}^{3}\partial^{l^{'}} \bf{E}.\partial_{l^{'}}\psi ) &   {\bf Id} 
	       {\partial_t \psi}  \\
               \end{array} 
 \right).
\end{equation}

\noindent Next, we compute the term $\psi \Dis\sum_{l=0}^{3}
 {\partial_l A^l}-2\psi S $, where $ S=1/2(C+C^{\;*})= C$. Note also that as $ A^l\;, 0\leq l\leq 3\;,$ are constant matrices, one 
has $\psi {\partial_l A^l}=0$.\par

\noindent Adding all the terms, we obtain

$$\{P,\psi\}+\sum^3_{l=0}\psi {\partial_l A^l}-2\psi S= 
\nonumber 
\left(
               \begin{array}{cccccccccccccccccccccc}   
               {\bf Id} {\partial_t \psi}  &   -
	       (\Dis\sum_{l'=1}^{3}\partial^{l^{'}} \bf{E}.\partial_{l^{'}}\psi )    
	       \\[0.3cm]\nonumber    
               (\Dis\sum_{l'=1}^{3}\partial^{l^{'}} \bf{E}.\partial_{l^{'}}\psi ) &   {\bf Id}
	       {\partial_t \psi}   \\\nonumber 
               \end{array} 
 \right)
-2\left(
               \begin{array}{cccccccccccccccccccccc}   
                  \psi{\bf Id}  &  0 \\[0.3cm] \nonumber        
                   0          &  0  \\ \nonumber 
                \end{array} 
 \right) \nonumber 
$$ 

\noindent and thus

\begin{equation}\label{3.21}   
\{P,\psi\}+\sum^3_{l=0}\psi {\partial_l A^l}-2\psi S =
\left(
               \begin{array}{cccccccccccccccccccccc}   
             ({\partial_t \psi}-2\psi ) {\bf Id}  & 
	     -(\Dis\sum_{l'=1}^{3}\partial^{l^{'}} \bf{E}.\partial_{l^{'}}\psi )    \\[0.3cm]        
              (\Dis\sum_{l'=1}^{3}\partial^{l^{'}} \bf{E}.\partial_{l^{'}}\psi ) &   
	      {\bf Id} {\partial_t \psi}   \\
               \end{array} 
 \right).
\end{equation}

\noindent Taking into account the form of $\scorpio$ deduced from the localisation property (\ref{localisation-constant}), we get

\begin{equation}\label{3.22}
\left\{
\begin{array}{ccc}
\displaystyle
<\scorpio,\{P,\psi\}+\sum^3_{l=0}\psi {\partial_l A^l}-2\psi S >=
< \left(
               \begin{array}{cccccccccccccccccccccc}   
                \zeta^{'}\otimes \zeta^{'}a(t,x,\zeta )    & 
		\zeta^{'}\otimes \zeta^{'} {c}(t,x,\zeta )  \\[0.3cm]        
              \zeta^{'}\otimes \zeta^{'}d(t,x,\zeta )      &  
	      \zeta^{'}\otimes \zeta^{'}b(t,x,\zeta )     \\
                \end{array} 
 \right),\\
\left(
               \begin{array}{cccccccccccccccccccccc}   
             ({\partial_t \psi}-2\psi ) {\bf Id}  
	     & -(\Dis\sum_{l'=1}^{3}\partial^{l^{'}} {\bf E}.\partial_{l^{'}}\psi )    \\[0.3cm]        
              (\Dis\sum_{l'=1}^{3}\partial^{l^{'}} {\bf E}.\partial_{l^{'}}\psi ) &   {\bf Id} 
	      {\partial_t\psi}   \\
               \end{array} 
 \right)>
\end{array}
\right.
\end{equation}

\noindent so that the propagation property reads as
\begin{equation}\label{ca1pro}
\left\{
\begin{array}{ccccccccccc}
\displaystyle
< \zeta^{'}\otimes \zeta^{'}a(t,x,\zeta ) ({\partial_t \psi}
-2\psi ) {\bf Id} +\zeta^{'}\otimes \zeta^{'} {c}(t,x,\zeta) 
(\Dis\sum_{l'=1}^{3} \partial^{l^{'}} {\bf E}.\partial_{l^{'}}\psi)>= 2 Re \mu_{{uf}_{11}}\;,\\[0.3cm]            
<-\zeta^{'}\otimes \zeta^{'} {a}(t,x,\zeta) (\Dis\sum_{l'=1}^{3} \partial^{l^{'}} {\bf E}.\partial_{l^{'}}\psi) 
+ \zeta^{'}\otimes \zeta^{'} {c}(t,x,\zeta) ){\bf Id} {\partial_t \psi}>=
2 Re \mu_{{uf}_{12}}\;,\\[0.3cm] 
<\zeta^{'}\otimes \zeta^{'}d(t,x,\zeta ) ({\partial_t \psi}
-2\psi ) {\bf Id} +
 \zeta^{'}\otimes \zeta^{'}b(t,x,\zeta ) (\Dis\sum_{l'=1}^{3} \partial^{l^{'}} {\bf E}.\partial_{l^{'}}\psi)>=
2 Re \mu_{{uf}_{21}}  \;, \\[0.3cm]    
<-\zeta^{'}\otimes \zeta^{'}d(t,x,\zeta ) (\Dis\sum_{l'=1}^{3} \partial^{l^{'}} {\bf E}.\partial_{l^{'}}\psi) +
\zeta^{'}\otimes \zeta^{'}b(t,x,\zeta ) {\bf Id} {\partial_t \psi}>=
2 Re \mu_{{uf}_{22}}
\end{array}
\right.
\end{equation}
with the notations explained in the statement of Theorem \ref{theoreme1}.\par
\noindent  Writing these equations in ${\cal D}'$ and then taking the trace of each equation, we get finally (\ref{ca1prooo}).\par
\noindent Now, it remains to take into account the propagation property coming from (\ref{systeme2-constant}). Since in this case, the characteristic polynomial is given by $B(\zeta ) =\Dis\sum_{\Dis j=1}^{3} B^{j}\zeta_{j}$, it follows with a small computation that one has the following propagation property
$$
\left\{
\begin{array}{ccc}
\displaystyle
< \left(
               \begin{array}{cccccccccccccccccccccc}   
                \zeta^{'}\otimes \zeta^{'}a(t,x,\zeta )    & 
		\zeta^{'}\otimes \zeta^{'} {c}(t,x,\zeta )  \\[0.3cm]        
              \zeta^{'}\otimes \zeta^{'}d(t,x,\zeta )      &  
	      \zeta^{'}\otimes \zeta^{'}b(t,x,\zeta )     \\
                \end{array} 
 \right),
\left(
               \begin{array}{cccccccccccccccccccccc}   
             \Gamma (\tilde x ,\zeta )  
	     & \Gamma (\tilde x,\zeta )\\[0.3cm]        
              \Gamma (\tilde x ,\zeta ) &  \Gamma (\tilde x ,\zeta) \\
               \end{array} 
 \right)>
\end{array}
\right. = 2Re \mu_{u\tilde\varrho}
$$
again using the notations in the statement of Theorem \ref{theoreme1}. Here $\Gamma (\tilde x,\zeta )$ is the $3\times 3$ matrix given by
$$\Gamma (\tilde x ,\zeta ) =\ diag (\partial_{x_1}\psi ,\partial_{x_2}\psi ,\partial_{x_3}\psi ).$$
\noindent Writing each equation and taking the trace, this gives the last constraint mentionned in the statement of Theorem \ref{theoreme1}.\par

\subsection{ Proof of Theorem \ref{theoreme2}: Non constant coefficient-scalar case}

\noindent Let us recall that we assume (\ref{hyp2}), that is $\ddot{\bf{\epsilon}}$, $\ddot{{\bf\eta}}$ and $\ddot{{\bf\sigma}} $ are $3\times 3$ scalar matrix valued functions given by

\begin{equation}\label{3.25}
\ddot{{\bf\epsilon}}={\bf \epsilon} ({\bf Id})_{3\times 3} \equiv \left(
               \begin{array}{cccccccccccccccccccccc}   
               \epsilon(x)&   0  &   0\\[0.3cm]     
                0&             \epsilon(x)  &   0\\[0.3cm]     
                0&             0  &   \epsilon(x)\\[0.3cm]     
                \end{array}                
 \right)
\end{equation} 

\noindent and

\begin{equation}\label{3.26}
\ddot{{\bf\eta}}={\bf\eta} ({\bf Id})_{3\times 3} \equiv  \left(
               \begin{array}{cccccccccccccccccccccc}   
               \eta (x)&    0  &    0\\[0.3cm]        
                0 &           \eta (x)   &   0\\[0.3cm]        
                0 &           0   &  \eta (x) \\[0.3cm]        
                \end{array}                
 \right)\;,
\ddot{{\bf\sigma}}={\bf\sigma} ({\bf Id})_{3\times 3} \equiv \left(
               \begin{array}{cccccccccccccccccccccc}   
               \sigma(x)&   0   &   0\\[0.3cm]        
               0&              \sigma(x)   &   0\\[0.3cm]        
               0&              0   &   \sigma(x)\\[0.3cm]        
                \end{array}                
 \right)
\end{equation}

\noindent where $ \epsilon$, $\eta$ and $\sigma $ are smooth functions in $C^{1}(\R^3) $, bounded from below.\par

\noindent The first equation of Maxwell's system can then again be written in the form of a symmetric 
system

\begin{equation}\label{systeme1-variable}
\Dis\sum_{i=0}^{3} A^i (x)\frac{\partial u^{\displaystyle\veps}}
{\partial x_{i}}+ C(x) u^{\displaystyle\veps}=
f^{\displaystyle\veps}
\end{equation}    

\noindent where

\begin{equation}\label{3.28}
 A^{\;0} =\left(
                  \begin{array}{cccccccccccccccccccccc}   
                     {\bf \epsilon}(x){\bf Id}  &   {\bf 0}  \\
                     {\bf 0}              &  {\bf \eta}(x){\bf Id}     \\

  \end{array}  
  \right)
\end{equation}

\noindent is a $(2\times 2)$ block matrix with $3\times 3$ blocks, the 
constant antisymmetric ${\bf Q_{\;k}}\;, 1\leq k\leq 3$ and the 
matrix $A^i\;,1\leq i \leq 3 $ being the same as in (\ref{elke}) and $C$ now given by 

\begin{equation}\label{3.29}
  C(x) =\left(
                  \begin{array}{cccccccccccccccccccccc}   
                     {\bf\sigma}(x) {\bf Id} &  {\bf 0}  \\
                     {\bf 0}           & {\bf 0}       \\

  \end{array}  
  \right)
\end{equation} 

\noindent is a $(2\times 2)$ block matrix, each block being of size $3\times 3$.\par
\noindent As the sequence $u^{\Dis\veps}$ converges weakly 
to zero in $L^2(\Omega )^6 $, again up to a subsequence, it defines an H- measure $\scorpio$, which is a $6\times 6$ matrix valued measure.\par


\noindent As in the preceding constant case, to express the localisation property linked with (\ref{systeme1-variable}), we compute the associated symbol $P(x,\zeta)$ which is here given by

\begin{equation}\label{3.31}
  P(x,\zeta )=\left(
                  \begin{array}{cccccccccccccccccccccc}   
                    \zeta_0 {\bf \epsilon}(x) {\bf Id}  & -{\bf E} \\[0.3cm]
                     {\bf E}             &  \zeta_0 {\bf\eta}(x) {\bf Id}  \\

  \end{array}  
  \right)
\end{equation}

\noindent with $ {\bf E}$ is still given by (\ref{matr}). The localisation property then states that $P \scorpio=0 $.\par
\noindent Let us first show the following
\noindent\begin{Lemma}\label{valeurs-propres}
\noindent Assume that $\zeta^{'}\neq 0$ and let $L(x,\zeta)=\Dis\sum_{j=1}^{3} A_{0}^{-1}(x) \zeta_j A^j $ 
be the dispersion matrix. Then $L$ has three eigenvalues, each with constant multiplicity two, given by 
$$\omega_{0}=0,\hspace{0.2cm} \omega_{+}=\zeta_{0}+v|\zeta^{'}|,\hspace{0.2cm} 
\omega_{-}=\zeta_{0}-v|\zeta^{'}|.$$ 
\noindent The  matrix $P'(x,\zeta)\equiv \zeta_{0} {\bf Id}+L(x,\zeta')$ has also 
three eigenvalues, given by 
$$\omega_{0}=\zeta_{0}\hspace{0.2cm}, \omega_{+}=\zeta_{0}+v|\zeta^{'}|\hspace{0.2cm}, 
\omega_{-}=\zeta_{0}-v|\zeta^{'}|,$$ 
\noindent each with constant multiplicity two, where $v(x)=\Dis\fa{1}{\sqrt{\epsilon(x)\eta(x)}}$ 
is the propagation speed. 
\end{Lemma}
\noindent {\bf Proof of Lemma \ref{valeurs-propres}}  

\noindent $L(x,\zeta^{'})$ can be rewritten as \par

\begin{equation}\label{definition-L}
L(x,\zeta^{'})=\Dis\sum_{j=1}^{3} A_{0}^{-1}(x) \zeta_{j} A^{j}=-
\left(
 \begin{array}{cccccccccccccccccccccc}   
0&0&0&0&-\zeta_{3}/\epsilon & \zeta_{2}/\epsilon  \\
0&0&0&\zeta_{3}/\epsilon&0&-\zeta_{1}/\epsilon    \\
0&0&0&-\zeta_{2}/\epsilon&\zeta_{1}/\epsilon&0\\
0&\zeta_{3}/\mu &-\zeta_{2}/\mu&0&0&0\\ 
-\zeta_{3}/\mu &0&\zeta_{1}/\mu &0&0&0\\  
\zeta_{2}/\mu &-\zeta_{1}/\mu&0&0&0&0\\    
\end{array}  
\right)
\end{equation}
\noindent or in block form 
\begin{equation}
L(x,\zeta^{'})=
\left(
 \begin{array}{cccccccccccccccccccccc}   
\bf{0}&-1/\epsilon \bf{E}\\
1/\mu \bf{E}&\bf{0}\\
\end{array}  
\right)\;.
\end{equation}
\noindent The action of the matrix $\bf{E}= {\bf E} (\zeta ')$ is also given as 
\begin{equation}
\bf{E}(\zeta^{'})[p]= \zeta^{'}\wedge p
\end{equation}
\noindent where $\bf{E}(\zeta^{'})$ is given in (\ref{matr}), for all $p\in \R^{3}$. Letting $\omega$ be an eigenvalue of $L(x,\zeta^{'})$ corresponding to an eigenvector 
$\bf{X}$, one has 
\begin{equation}
\left(
 \begin{array}{cccccccccccccccccccccc}   
\bf{0}& -1/\epsilon \bf{E}\\
1/\mu \bf{E}&\bf{0}\\
\end{array}  
\right){\bf X}=\omega {\bf X}
\end{equation}
\noindent or 
\begin{equation}\label{vp1}
\left(
 \begin{array}{cccccccccccccccccccccc}   
\bf{0}& -1/\epsilon \bf{E}\\
1/\mu \bf{E}&\bf{0}\\
\end{array}  
\right)
\left(
 \begin{array}{cccccccccccccccccccccc}   
\bf{u}\\
\bf{v}\\
\end{array}  
\right)=
\omega \left(
 \begin{array}{cccccccccccccccccccccc}   
\bf{u}\\
\bf{v}\\
\end{array}  
\right).
\end{equation}
\noindent Let us check that $\omega=0$ is an eigenvalue. In order to see this, we have to solve the algebric system
\begin{equation}
\left\{
\begin{array}{ccccccc}
-1/\epsilon \bf{E}(\zeta^{'})\vec{v}=0,\\[0.3cm]
1/\mu \bf{E}(\zeta^{'})\vec{u}=0 
\end{array}  
\right.
\end{equation}
\noindent which is equivalent to
\begin{equation}
\left\{
\begin{array}{ccccccc}
\vec{\zeta^{'}}\wedge \vec{\bf{v}}=0,\\[0.3cm]
\vec{\zeta^{'}}\wedge \vec{\bf{u}}=0.\\
\end{array}  
\right.
\end{equation}

\noindent It follows that $\vec {\bf u}$ and $\vec {\bf v}$ are colinear to $\zeta '$, and thus that $(\vec 0 , \vec \zeta ')$ and $(\vec \zeta ',\vec 0 )$ 
is a basis of the eigenspace corresponding to the eigenvalue $0$ of $L$. It shows also that $0$ is indeed an eigenvalue of multiplicity two.\par
\noindent It remains to find the eigenvalues $\omega \neq 0$ of $L$. For this purpose, set
$$D(\zeta^{\;'}):= 
\Dis\sum_{j=1}^{3} \frac{\zeta_{\;j}}{| \zeta^{\;'}|} A^{j} =(\Dis\fa{1}{|\zeta|} L(x,\zeta^{'})).$$ 
\noindent We can as well assume that $\zeta^{\;'}\in S^{\;2}$, the unit sphere in $\R^{3}$. Then $\omega =\mp 1$ are 
eigenvalues of the matrix $ D(\zeta^{\;'}) $.\par 

\noindent Indeed, we have to solve the system of equation
\begin{equation}
\left\{
\begin{array}{ccccccc}
- \vec{\zeta^{'}}\wedge \vec{\bf{v}} =\omega \vec{\bf{u}}\;,\\[0.3cm]
\vec{\zeta^{'}}\wedge \vec{\bf{u}}=\omega \vec{\bf{v}}\;.\\
\end{array}  
\right.
\end{equation}

\noindent For $\omega=1$ to be an eigenvalue of the matrix $D(\zeta^{\;'})$, we have to solve the system

\begin{equation}\label{3.42}
\left\{
\begin{array}{ccc}
\displaystyle 
-\vec{\zeta^{'}} \wedge  \vec{\bf{v}} = \vec{\bf{u}}\;,\\[0.3cm]
\vec{\zeta^{'}} \wedge  \vec{\bf{u}}= \vec{\bf{v}}
\end{array}
\right.
\end{equation}
\noindent which admits a $2D$ space of solutions, and for $\omega =-1$ to be an eigenvalue of the matrix $D(\zeta^{\;'})$, we have to solve

\begin{equation}\label{3.43}
\left\{
\begin{array}{ccc}
\displaystyle 
-\vec{\zeta^{'}}\wedge  \vec{\bf{v}} =- \vec{\bf{u}}\;,\\[0.3cm]
\vec{\zeta^{'}}\wedge  \vec{\bf{u}}=-\vec{\bf{v}}\;
\end{array}
\right.
\end{equation}
\noindent which again admits a $2D$ space of solutions. Thus, all in all, we have obtained that
$$\omega_{\;1} =-1\;,\omega_{\;3} =1 $$
\noindent are the eigenvalues of multiplicity $m=2 $ of the matrix $ D(\zeta^{\;'})$. Recall that $ \omega_{\;2} =0$ is 
also an eigenvalue of the matrix $ D(\zeta^{\;'})$ with multiplicity $m=2 $.\par

\noindent Next, note that if $\lambda $ is an eigenvalue of a matrix $A$ 
then $c\lambda $ is an eigenvalue of the matrix $cA$, 
where $c>0$ is a constant . Thus we conclude that
\begin{equation}\label{val}
- |\zeta^{\;'}|\;, 0\;, |\zeta^{\;'}|\; 
\end{equation}
\noindent are the eigenvalues of the matrix $B(\zeta^{\;'})=L(x,\zeta^{'})$ with multiplicity $m=2$,  for all 
$\zeta^{\;'}\neq 0 $.\par 
\noindent Set $ \mu(x,\zeta ) =\lambda(x,\zeta^{'})-\zeta_0 $ and 
recall that if $\lambda(x,\zeta) $ is an eigenvalue of the matrix 
$P$ corresponding to an eigenvector $ \vec{u}=\vec{u}(x,\zeta )$, 
then $\mu(x,\zeta ) =\lambda(x,\zeta^{'})-\zeta_0 $ is an eigenvalue 
of the matrix $L(x,\zeta^{'})$ corresponding to an eigenvector 
$ \vec{\bf{X}}=\vec{\bf{X}}(x,\zeta^{\;'})$. If 
$\mu(x,\zeta ) =\lambda(x,\zeta )-\zeta_0 $ is an 
eigenvalue of the matrix  $L(x,\zeta^{'})$  corrsponding to an 
eigenvector $ \vec{\bf{X}}=\vec{\bf{X}}(x,\zeta^{\;'})$, then 
$\lambda(x,\zeta ) =\mu(x,\zeta ) +\zeta_0 $ is an 
eigenvalue of the matrix $P'(x,\zeta )$ corresponding to an eigenvector 
$ \vec{\bf{X}}=\vec{\bf{X}}(x,\zeta)$. As the eigenvalues of the matrix 
$L(x,\zeta^{'})$ are given in (\ref{val}), we can conclude that
\begin{equation}\label{val2}
\zeta_{\;0}-v|\zeta^{\;'}|\;, \zeta_{\;0}\;, 
\zeta_{\;0}+v|\zeta^{\;'}|  
\end{equation}
are the eigenvalues of the matrix $P'(x,\zeta )$ with 
multiplicity $m=2 $. Now we shall show that the 
propagation speed $v$ is given by 
$$v(x)=\Dis\fa{1}{\sqrt{\epsilon(x)\mu(x)}}\;.$$
\noindent Indeed, from (\ref{vp1}) and for $\omega\neq 0$, one has 
\begin{equation}
\left\{
\begin{array}{ccccccc}
-1/\epsilon \bf{E}(\zeta^{'})\vec{v}=\omega\vec{u} \\
1/\mu \bf{E}(\zeta^{'})\vec{u}=\omega\vec{v} 
\end{array}  
\right.
\end{equation}
\noindent or
\begin{equation}
\left\{
\begin{array}{ccccccc}
\vec{\zeta^{'}}\wedge \vec{\bf{v}}=-\epsilon \omega \vec{\bf{u}}\\
\vec{\zeta^{'}}\wedge \vec{\bf{u}}=\mu \omega \vec{\bf{v}}\\
\end{array}  
\right.
\end{equation}
\noindent and thus 
$$\vec{\bf{v}}=\Dis\fa{1}{\omega \mu }\vec{\zeta^{'}}\wedge \vec{\bf{u}}\;.$$
\noindent Thus
\begin{equation}
\left\{
\begin{array}{ccccccc}
\vec{\zeta^{'}}\wedge (\fa{1}{\mu \omega } \vec{\zeta^{'}}\wedge \vec{\bf{u}})=
-\epsilon \omega\vec{\bf{u}}\Rightarrow \vspace{0.2cm}\\
\vec{\zeta^{'}}\wedge (\vec{\zeta^{'}}\wedge \vec{\bf{u}})= 
-\epsilon \mu \omega^{2} \vec{\bf{u}}\Rightarrow \vspace{0.2cm} \\
\{ \vec{\zeta^{'}}\wedge (\vec{\zeta^{'}}\wedge \vec{\bf{u}}) \}.\vec{\bf{u}}=
-\epsilon \mu \omega^{2}||\vec{\bf{u}}||^{2}\Rightarrow\vspace{0.2cm} \\
(\vec{\zeta^{'}}\wedge \vec{\bf{u}}).(\vec{\zeta^{'}}\wedge \vec{\bf{u}})= 
-\epsilon \mu \omega^{2}||\vec{\bf{u}}||^{2}\Rightarrow \vspace{0.2cm} \\
||\vec{\zeta^{'}}\wedge \vec{\bf{u}}||^{2}=- \epsilon \mu \omega^{2} ||\vec{\bf{u}}||^{2}\\
\end{array}  
\right.
\end{equation}
\noindent but $\vec{\zeta^{'}}$ is orthogonal on the vector $\vec{\bf{u}}$ and thus 
\begin{equation}
\omega^{2}=\Dis\fa{||\zeta^{'}||^{2} ||\vec{\bf{u}}||^{2}}{\epsilon \mu ||\vec{\bf{u}}||^{2}}\Longrightarrow 
\omega=\pm \Dis\fa{||\zeta^{'}||}{\sqrt{\epsilon \mu }}=\pm \Dis\fa{|\zeta^{'}|}{\sqrt{\epsilon \mu }}=
v|\zeta^{'}|\;.
\end{equation}
\noindent This ends the proof of Lemma \ref{valeurs-propres}.\par
\smallskip
\noindent For the eigenvectors of the matrix $P'(x,\zeta )$, one first chooses an orthonormal basis 
of $\R^{3}$. This basis consists of the propagation triple in the direction of 
propagation $\hat{\zeta '}$ and of the two transverse unit vectors $ z^{(1)}(\zeta '),z^{(2)}(\zeta ')$. \par 
\noindent Let $(\hat{\zeta '},z^{(1)}(\zeta '),z^{(2)}(\zeta ')) \in \R^{3}$ be this basis of propagation. 
In polar coordinates they are given by, see for more details \cite{Geoleon}, \cite{Geoleonkel}
\begin{equation}\label{vecteurs-propres1}
\hat{\zeta '}= {{\zeta '}\over{\mid\zeta '\mid}} =
\left(
 \begin{array}{cccccccccccccccccccccc}   
\sin \theta \cos \phi \\
\sin \theta \sin \phi \\
\cos \theta 
\end{array}  
\right)\;,\hspace{0.2cm}
z^{(1)}(\zeta ')=
\left(
 \begin{array}{cccccccccccccccccccccc}   
\cos\theta \cos \phi \\
\cos \theta \sin \phi \\
-\sin \theta \\
\end{array}  
\right)\;,\hspace{0.2cm}
z^{(2)}(\zeta ')=
\left(
 \begin{array}{cccccccccccccccccccccc}   
-\sin \phi \\
\cos \phi \\
0\\ 
\end{array}  
\right)\;
\end{equation}

\noindent where $|\zeta '|=(\zeta_{1}^{2}+\zeta_{2}^{2}+\zeta_{3}^{2})^{1/2}$.\par

\noindent Then, one can show that the eigenvectors of the matrix $P'(x,\zeta )$ are given by 
\begin{equation}\label{vecteurs-propres2}
\left\{
\begin{array}{ccccccc}
b^{1}_0=\Dis\fa{1}{\sqrt \epsilon }(\hat{\zeta '},0) \;,\hspace{0.3cm} b^{2}_0=\Dis\fa{1}{\sqrt{\mu}}(0,\hat{\zeta '}) \vspace{0.2cm}\\
b_{+}^{1}=(\Dis\fa{1}{\sqrt{2\epsilon}} z^{1},\Dis\fa{1}{\sqrt{2\mu }}z^{2})\;,\hspace{0.3cm}
b_{+}^{2}=(\Dis\fa{1}{\sqrt{2\epsilon}} z^{2},-\Dis\fa{1}{\sqrt{2\mu}}z^{1})\vspace{0.2cm} \\   
b_{-}^{1}=(\Dis\fa{1}{\sqrt{2\epsilon}} z^{1},-\Dis\fa{1}{\sqrt{2\mu }}z^{2})\;,\hspace{0.3cm}
b_{-}^{2}=(\Dis\fa{1}{\sqrt{2\epsilon}} z^{2},\Dis\fa{1}{\sqrt{2\mu}}z^{1})\\   
\end{array}  
\right.
\end{equation}
\noindent The eigenvectors $b^{1}_0$ and $b^{2}_0$ represent the non-propagating longitudinal and the other eigenvectors 
correspond to transverse modes of propagation with respect the propagation speed $v$.\par
\smallskip
\noindent\begin{Lemma}\label{decomposition-variable} The H-measure $ \scorpio $ has the form
\begin{equation}\label{rhande} 
\scorpio= b^{1}_0\otimes b^{1}_0 a_{0}+ b^{2}_0\otimes  b^{2}_0 b_{0}+
b_{+}^{1}\otimes   b_{+}^{1} a_{+}+ b_{+}^{2}\otimes  b_{+}^{2} b_{+}+
b_{-}^{1}\otimes  b_{-}^{1} a_{-}+ b_{-}^{2}\otimes  b_{-}^{2} b_{-}
\end{equation}
\smallskip 
\noindent where $a_{0}$ ,$ b_{0}$ are two positives measures supported in 
the set $\{\zeta_0 =0 \}$, $ a_{+}$, $b_{+}$ are two positives measures supported in the set $\{\zeta_0 =-v\mid\zeta '\mid\}$, and $ a_{-}$, $b_{-}$ are two positives measures 
supported in the set $\{\zeta_0 =+v\mid\zeta '\mid\}$. Here $b^{1}_0$, $b^{2}_0$, $b_{+}^{1}$, $b_{+}^{2}$, $b_{-}^{1}$, and $b_{-}^{2}$ are the 
eigenvectors of the matrix $P'$ given by (\ref{vecteurs-propres1}) and (\ref{vecteurs-propres2}).\par

\end{Lemma}   
\noindent{\bf Proof of Lemma (\ref{decomposition-variable})}         
\noindent Since $P\scorpio =$ by the localization property, it follows that
\begin{equation}\label{rof}
(A_{0})^{-1}P \scorpio =0\;.
\end{equation}
\noindent Since $P'(x,\zeta )=(A_{0})^{-1}P(x,\zeta )$, thus (\ref{rof}) becomes 
\begin{equation}\label{roof}
P'(x,\zeta ) \scorpio =0\;.
\end{equation}
\noindent It follows that the support of the H-measure $\scorpio $
is included in the set 
\begin{equation}\label{roffo}
U=\{(x,\zeta)\in \R^{4}\times S^{3}\;, \mbox{det} P'=0\}\;. 
\end{equation}

\noindent Note that, for every $t=\zeta_{0}\;, x\;,\zeta $ fixed, the matrix  
$ P'$ is diagonalizable. In fact, we will discuss the following cases :\par 

\noindent {\bf a) Case of equal eigenvalues} \par
\noindent i) If $\omega_{\;0}=\omega_{\;+}$, (resp.$\omega_{\;0}=\omega_{\;-}$), then from Lemma \ref{valeurs-propres}), 
one has $\zeta'=0$, (resp.$\omega_{\;0}=\omega_{\;-}$), and thus $\zeta_{\;0}=\pm 1$, 
and $P'=\zeta_{\;0}Id$, which is diagonal.\par 
\noindent ii) If $\omega_{\;0}\neq\omega_{\;+}$, $\omega_{\;0}\neq\omega_{\;-}$, but 
$\omega_{\;+}=\omega_{\;-}$. In this case, again from Lemma \ref{valeurs-propres}, one has 
 $\zeta'=0$, and thus $\zeta_{\;0}=\pm 1$ and $P'=\zeta_{\;0}Id$, which is diagonal.\par
\noindent {\bf b) Case of distinct eigenvalue} \par
\noindent If $\omega_{\;0}\neq\omega_{\;+}\neq \omega_{\;-}$, then again from Lemma 
(\ref{valeurs-propres}), using the basis of eigenvectors corresponding to an 
the eigenvalues $\omega_{\;0}$, $\omega_{\;+}$, $\omega_{\;-}$, given by (\ref{vecteurs-propres1}) and (\ref{vecteurs-propres2}), $P'=\zeta_{\;0}Id$ is diagonal.\par 
\noindent Next, as the support of the H-measure $\scorpio $ is contained in the set of points
(\ref{roffo}), and as the matrix $P'$ is diagonalizable, recalling that the determinant of the matrix $P'$ 
is given by  (eventually with powers)
$$\mbox{det} P'\sim \omega_{0} \omega_{+} \omega_{-}\;.$$  
It follows that the support of the H-measure $\scorpio $ is included in the set 
$$U=\{(x,\zeta)\in \R^{4}\times S^{3}\;, \omega_{\;0}=0\cup \omega_{+}=0\cup \omega_{-}\}\;. $$
\noindent Next, one has
$$P'\scorpio_{j}=0\hspace{0.6cm}  j=1,2,.....6$$
\noindent where  $\scorpio_j$ is the $j-th$ column vector of the matrix $\scorpio$. Using Lemma (\ref{valeurs-propres}), we get that, for $\omega_{0}=0$, there exist two scalars 
$\alpha_{j}\;\;, \beta_{j}$ such that 
$$\scorpio_{j}= \alpha_{j} b^{1}_0+\beta_{j} b^{2}_0\;.$$ 
\noindent Repeating the sames steps for the other eigenvalues $\omega_{\;+}\;, \omega_{\;-}$, and using the hermitian property of the H-measure, this {\bf ends the proof of Lemma \ref{decomposition-variable}}.\par
\smallskip   
\noindent In order to write down the propagation property for equation (\ref{systeme1-variable}), we 
first need to compute the Poisson bracket, which in this case, is given by

\begin{equation}\label{3.44}
\{ P,\psi\} =\left(
               \begin{array}{cccccccccccccccccccccc}   
             {\bf\epsilon}(x){\bf Id} {\partial_t \psi}-\zeta_0 
\Dis\sum_{l'=1}^{3}\partial^{\;l^{'}} \psi \partial_{\;l^{'}}{\bf\epsilon}(x){\bf Id} & -(\Dis\sum_{l'=1}^{3}\partial^{\;l^{'}} 
{\bf E}.\partial_{l^{\;'}}\psi )    \\[0.3cm]        
              (\Dis\sum_{l'=1}^{3}\partial^{\;l^{'}} {\bf E}.\partial_{l^{\;'}}\psi ) &  {\bf\eta}(x){\bf Id} 
{\partial_t \psi}-\zeta_0 \Dis\sum_{l'=1}^{3}\partial^{\;l^{'}} \psi 
\partial_{\;l^{'}}{\bf\eta}(x){\bf Id}   \\
               \end{array} 
 \right)\;.
\end{equation}


\noindent It follows that, since $\Dis\sum_{k=0}^{3}\partial_{k} A^{k}=0$


\begin{equation}\label{3.46}
\left\{
\begin{array}{cccccc}
\Dis
\{P,\psi\}+\psi \sum^3_{k=0}{\partial_k A^k}-2\psi S = \\[0.3cm]
= \left(
               \begin{array}{cccccccccccccccccccccc}   
               {\bf\epsilon}(x){\bf Id}{\partial_t \psi}-\zeta_0 
\Dis\sum_{l'=1}^{3}\partial^{\;l^{'}} \psi \partial_{\;l^{'}}{\bf\epsilon}(x) {\bf Id}-2\psi\sigma {\bf Id}  & 
-(\Dis\sum_{l'=1}^{3}\partial^{\;l^{'}} {\bf E}.\partial_{l^{\;'}}\psi )    \\[0.3cm]        
              (\Dis\sum_{l'=1}^{3}\partial^{\;l^{'}} {\bf E}.\partial_{l^{\;'}}\psi ) &    {\bf\eta}(x){\bf Id}
{\partial_t \psi}-\zeta_0 \Dis\sum_{l'=1}^{3}\partial^{\;l^{'}} \psi 
\partial_{\;l^{'}}{\bf\eta}(x){\bf Id}     \\
               \end{array} 
 \right)\;.
\end{array}
\right.
\end{equation}
\noindent Writing the H-measure $\scorpio$ as
\begin{equation}\label{forme-bloc-variable}
\scorpio = \left(\begin{array}{cccccccccc} \scorpio_{11} & \scorpio_{12} \\ \scorpio_{21} &\scorpio_{22}\\ \end{array}\right)
\end{equation}
where $\scorpio_{ij}$ are $3\times 3$ matrix valued measures, it follows that one has
\begin{equation}\label{propagation-variable}
\left\{\matrix{\displaystyle 
-\Dis\veps (x) \partial_t \scorpio_{11} +\zeta_0\sum^3_{l'=1} \partial_{l'}\veps (x) \partial^{l'}\scorpio_{11} - 2\sigma \scorpio_{11} - \sum^{3}_{l'=1} \partial^{l'} {\bf E} .\partial_{l'} \scorpio_{12} = 2 Re \mu_{{uf}_{11}},\cr
\displaystyle -\eta (x) \scorpio_{12} +\zeta_0 \sum^3_{l'=1} \partial_{l'} \eta (x) \partial^{l'} \scorpio_{12} +\sum^3_{l' =1} \partial^{l'} {\bf E} \partial_{l'}\scorpio_{11} = 2 Re \mu_{{uf}_{12}}, \cr
\displaystyle 
-\Dis\veps (x) \partial_t \scorpio_{21} +\zeta_0\sum^3_{l'=1} \partial_{l'}\veps (x) \partial^{l'}\scorpio_{11} - 2\sigma \scorpio_{21} - \sum^{3}_{l'=1} \partial^{l'} {\bf E} .\partial_{l'} \scorpio_{22} = 2 Re \mu_{{uf}_{21}},\cr
\displaystyle -\eta (x) \scorpio_{22} +\zeta_0 \sum^3_{l'=1} \partial_{l'} \eta (x) \partial^{l'} \scorpio_{22} +\sum^3_{l' =1} \partial^{l'} {\bf E} \partial_{l'}\scorpio_{21} = 2 Re \mu_{{uf}_{22}} \\
}\right.
\end{equation}
which is the form given in Theorem \ref{theoreme2}.\par
\noindent Using the eigenvector basis (\ref{vecteurs-propres1}) and (\ref{vecteurs-propres2}), with the decomposition given by Lemma \ref{decomposition-variable}, we note that the elements of the H-measure $\scorpio$ defined in (\ref{forme-bloc-variable}), can be expressed as
\begin{equation}\label{decomposition-bloc}
\left\{
\begin{array}{cccccc}
\displaystyle
\scorpio_{11}=
       \fa{1}{\epsilon} [(\hat{\zeta '}\otimes \hat{\zeta '})a_{0}+
\fa{1}{2}(z^{1}\otimes z^{1}) a_{+}+\fa{1}{2}(z^{2}\otimes z^{2}) b_{+}
+\fa{1}{2}(z^{1}\otimes z^{1})a_{-}\\[0.2cm]
+\fa{1}{2}(z^{2}\otimes z^{2}) b_{-}],\\[0.3cm]
\scorpio_{12}=\fa{v}{2}
[(z^{1}\otimes z^{2})a_{+}-
(z^{2}\otimes z^{1}) b_{+}-(z^{1}\otimes z^{2})a_{-}+
(z^{2}\otimes z^{1}) b_{-}],\\[0.3cm]
\scorpio_{21}=
\fa{v}{2}[(z^{2}\otimes z^{1}) a_{+}-(z^{1}\otimes z^{2})b_{+}-
(z^{2}\otimes z^{1})a_{-}+(z^{1}\otimes z^{2})b_{+}],\\[0.3cm]
\scorpio_{22}= 
\fa{1}{\mu}[(\hat{\zeta '}\otimes \hat{\zeta '})b_{0}+\fa{1}{2}(z^{2}\otimes z^{2}) a_{+}+\fa{1}{2}(z^{1}\otimes z^{1})b_{+}
\\[0.2cm]+
\fa{1}{2}(z^{2}\otimes z^{2})a_{-}+\fa{1}{2}(z^{1}\otimes z^{1})b_{-} ].
\end{array}
\right.\;
\end{equation}
If one wants to find an equation for $a_0$ for instance, one can proceed as follows. Recalling the form of $\scorpio_{11}$ just given above, we take the equation for it in (\ref{propagation-variable}) and apply it to the vector $\zeta '$.\par
\noindent Finally, because the divergence constraint in Maxwell's system only involves scalar valued functions, similar statement as in the end of Theorem \ref{theoreme1} holds true again.\par
\maketitle 
\noindent { \bf Acknowledgements: } The author would like to thank Radjesvarane Alexandre for several
discussions and suggestions during the preparation of this paper. \par

\end{document}